\documentclass[10pt]{article}

\usepackage{amsfonts,amsmath,amssymb,colortbl}

\begin{document}
\title{Enumerations of maximum partial triple systems on 16 and 17 points}
\author{Fatih Demirkale\\
Department of Mathematics\\
Y{\i}ld{\i}z Technical University\\
Esenler, \.{I}stanbul 34220\\ TURKEY\\
\texttt{(fatihd@yildiz.edu.tr)}\vspace{-1mm}
\and
Diane Donovan\\
Centre for Discrete Mathematics and Computing\\
University of Queensland, St. Lucia 4072\\ AUSTRALIA\\
\texttt{(dmd@maths.uq.edu.au)}\vspace{-1mm}
\and
Mike Grannell\footnote{Corresponding author}\\
School of Mathematics and Statistics \\
The Open University, Walton Hall\\
Milton Keynes MK7 6AA\\
UNITED KINGDOM\\ \texttt{(m.j.grannell@open.ac.uk)}}

\maketitle
\vspace{-8mm}
\begin{abstract}
For $v\equiv 1$ or 3 (mod 6), maximum partial triple systems on $v$ points are Steiner triple systems, STS($v$)s. The 
80 non-isomorphic STS(15)s were first enumerated around 100 years ago, but the next case for Steiner triple systems was 
unresolved until around 2004 when it was established that there are precisely 11\,084\,874\,829 non-isomorphic 
STS(19)s. In this paper we complete enumeration of non-isomorphic maximum partial triple systems for $v\le 19$. It is 
shown that there are 35\,810\,097 systems on 17 points and 47\,744\,568 on 16 points. We also establish that there are 
precisely 157\,151 non-isomorphic pairwise balanced designs, PBD($17,\{3,5\}$)s, having a single block of size 5. 
Structural properties of these systems are determined, including their automorphism groups, and the numbers of Pasch 
configurations, mitres and Fano planes contained in them. The systems themselves are available from the authors.
\end{abstract}

\section{Introduction}\label{Intro}
We will be concerned with certain types of combinatorial designs, in particular partial triple systems. A \emph{partial 
triple system} of order $v$ and index $\lambda$, PTS($v,\lambda$), is a pair $(V,\mathcal{B})$, where $V$ is a set of 
$v$ elements (the \emph{points}) and $\mathcal{B}$ is a collection of 3-element subsets of $V$ (the \emph{triples} or 
\emph{blocks}) with the property that every 2-element subset of $V$ occurs at most $\lambda$ times amongst the triples 
of $\mathcal{B}$. A PTS($v,\lambda)=(V,\mathcal{B})$ is said to be \emph{maximal} if there is no other 
PTS($v,\lambda)=(V,\mathcal{B}')$ with $\mathcal{B}\subseteq\mathcal{B}'$. Informally, a 
PTS($v,\lambda)=(V,\mathcal{B})$ is maximal if no further triples can be added to $\mathcal{B}$ without violating the 
condition that pairs appear at most $\lambda$ times. A maximal PTS($v,\lambda$) is denoted as MPTS($v,\lambda$). A 
\emph{maximum} MPTS($v,\lambda$), denoted as MMPTS($v,\lambda$), is an MPTS($v,\lambda$) having the greatest number of 
triples of any MPTS($v,\lambda$). Rather than calling an MMPTS($v,\lambda$) a ``maximum maximal partial triple system'', 
we prefer to use the term ``maximum partial triple system''. In the case $\lambda=1$, we will denote these systems as 
PTS($v$), MPTS($v$) and MMPTS($v$) respectively, and note that $\mathcal{B}$ is then necessarily a \emph{set} of 
triples, i.e. no triple can be repeated. The \emph{leave} of an PTS($v$) (or MPTS($v$), MMPTS($v$)) is the set of pairs 
of points that do not appear in any triple of the system.

If $K$ is a set of positive integers, then a \emph{pairwise balanced design}, \linebreak PBD($v,K,\lambda)$, of order 
$v$ and index $\lambda$, is a pair $(V,\mathcal{B})$, where $V$ is a set of $v$ elements (the \emph{points}) and 
$\mathcal{B}$ is a collection of subsets of $V$ (the \emph{blocks}) with the properties that every 2-element subset of 
$V$ occurs precisely $\lambda$ times amongst the blocks of $\mathcal{B}$, and if $B\in\mathcal{B}$ then $|B|\in K$. The 
elements of $K$ are called the block sizes. If $\lambda=1$, such a system is denoted more simply as PBD($v,K$). If the 
leave of a PTS($v$) is combined with its set of triples, then a pairwise balanced design PBD($v,\{2,3\}$) is obtained. 
Thus there is a rough equivalence between pairwise balanced designs with block sizes 2 and 3, and partial triple 
systems.

A \emph{Steiner triple system} of order $v$, STS($v$), is an MMPTS($v$) in which each pair of points appears in some 
triple, that is to say the leave is empty. It is well known that the necessary and sufficient condition for the 
existence of an STS($v$) is that $v\equiv 1$ or 3 (mod 6), and these values are called \emph{admissible}. When 
$v\not\equiv 1$ or 3 (mod 6), an STS($v$) does not exist and an MMPTS($v$) is, in a sense, the closest we can get to a 
Steiner triple system. When $v\equiv 5$ (mod 6) an MMPTS($v$) necessarily has a leave comprising four pairs of the form 
$\{a,b\},\{b,c\},\{c,d\},\{d,a\}$ (see \cite[page 553]{HBK}). In this case, taking the pairs as edges of a graph on the 
four points, the leave can be represented as a 4-cycle $(a,b,c,d)$. For the cases $v\equiv 0$ or 2 (mod 6), an 
MMPTS($v$) corresponds to an STS($v+1$) in which one point $a$ has been deleted, resulting in a leave comprising $v/2$ 
disjoint pairs. For the case $v\equiv 4$ (mod 6) a similar correspondence holds between an MMPTS($v$) and an 
MMPTS($v+1$) where a point $a$ from the leave of the latter has been deleted giving the former a leave consisting of the 
three intersecting pairs $\{b,c\},\{c,d\},\{c,e\}$ and, disjoint from these, a further $(v-4)/2$ disjoint 
pairs. 

It should be noted that maximum partial triple systems, MMPTS($v$)s, are also called \emph{optimal} $2-(v,3,1)$ 
\emph{packings} (see \cite[Section VI.40]{HBK}). For $v\equiv 0,1,2,3$ (mod 6), they are additionally known as 
$(v_r,b_3)$ \emph{configurations}, where $b$ is the number of triples and $r$ is the number of triples containing each 
point (see \cite[Section VI.7]{HBK}).

Two designs, $(V_1,\mathcal{B}_1)$ and $(V_2,\mathcal{B}_2)$ (of the same order $v$) are said to be isomorphic if 
there is a mapping $\phi$ from $V_1$ onto $V_2$ that takes blocks of $\mathcal{B}_1$ onto those of $\mathcal{B}_2$. 
Such a mapping $\phi$ is called an isomorphism. An isomorphism from a design to itself is called an automorphism, and 
the set of all automorphisms of a given design forms a group under composition. This group is referred to as the
automorphism group of the design. In discussing enumerations, it is generally the number of isomorphism classes of each 
design that is sought, in which case the point set may be fixed, for example as the set of integers 
$\{0,1,\ldots,v-1\}$.

As regards enumeration of these designs for small orders $v$, the situation for Steiner triple systems (up to 
isomorphism) is as follows. There are unique STS($v$)s for $v=3$, $v=7$ and $v=9$. There are two STS(13)s and 80 
STS(15)s; the latter result \cite{WCC} dates from 1919 and it was certainly not trivial to obtain this in the days 
prior to the advent of electronic computers. It took a further 85 years before the STS(19)s were enumerated in a paper 
by Kaski and \"{O}sterg{\aa}rd \cite{KO}; there are 11\,084\,874\,829 such systems and this result was obtained by a 
computer search. It seems unlikely that the next case, STS(21), will be enumerated any time soon since it is estimated 
that the number of these is several orders of magnitude greater than in the STS(19) case. It is known that the number 
of STS($v$)s is $v^{v^2(\frac16-o(1))}$ as $v\to \infty$ \cite{Wilson} (see also \cite{CR} Section 5.2.), a result 
that extends to MMPTS($v$)s. 

In the case when $v\equiv 5$ (mod 6), the MMPTS(5) is unique up to isomorphism, but there are two MMPTS(11)s 
\cite{CR1} and these are given in Table \ref{Table 1}. There are generally two types of MMPTS($v$) designs for $v\equiv 
5$ (mod 6). The first type is when the leave $\{a,b\},\{b,c\},\{c,d\},\{d,a\}$ is accompanied by a pair of intersecting 
triples of the design, $\{a,c,e\}$ and $\{b,d,e\}$. Putting the 6 pairs from these two triples together with the leave 
gives a copy of the complete graph $K_5$, and we will refer to this as the \emph{quintuple} case, and such designs as 
being of Type-Q. The second type is when the leave is accompanied by a pair of disjoint triples of the design, 
$\{a,c,e\}$ and $\{b,d,f\}$ ($e\ne f$), and we will refer to this as the \emph{non-quintuple} case, and such designs 
as being of Type-N. In Table \ref{Table 1}, the design $Q$ is of Type-Q, and the design $N$ is of Type-N. In 
both cases the point set is $\{0,1,\ldots,10\}$ (with 10 denoted by $A$) and the leave is the 4-cycle $(0,1,2,3)$. Here 
and elsewhere, when listing pairs and triples, set brackets and/or commas may be omitted, so that $\{a,b,c\}$ may be 
written as $abc$. A reader with some patience can construct these two designs by hand and verify that there is 
precisely one of each type up to isomorphism. 
\begin{table}[htbp]
\begin{center}
\[\begin{array}{|c|l|}\hline \mbox{System} & \mbox{triples}\\ \hline
Q & 024, 056, 078, 09A, 134, 159, 168, 17A, 257, 269, 28A,\\ &  35A, 367, 389, 458, 46A, 479.\\
\hline
N & 024, 057, 068, 09A, 135, 148, 169, 17A, 259, 26A, 278,\\ &  34A, 367, 389, 456, 479, 58A.\\ 
\hline
\end{array}\]
\caption{The two maximum partial triple systems of order 11.\label{Table 1}}
\end{center}
\end{table}

Given the correspondence between MMPTS($v$)s and MMPTS($v+1$)s, we begin by enumerating MMPTS(17)s. Our results then 
enable us to enumerate the MMPTS(16)s. We first show that there are precisely 35\,810\,097 MMPTS(17)s. We find that 
there are 2\,350\,733 MMPTS(17)s of Type-Q, and 33\,459\,364 \linebreak MMPTS(17)s of Type-N. These results were 
obtained by a computer search described in Section \ref{Searching17}. The next value of $v\equiv 5$ (mod 6) is $v=23$. 
The number of MMPTS(23)s is likely to be further orders of magnitude greater than the number of STS(21)s and so will 
probably remain out of reach for a considerable period of time.

There is a general construction of MMPTS($v$)s for $v\equiv 5$ (mod 6) from PBD($v,\{3,5\}$)s having precisely one 
block of size 5; such designs (denoted as PBD($v,\{3,5^*\}$)s) exist for all $v\equiv 5$ (mod 6) (see \cite{CR} Section 
6). If the block of size 5 is $\{a,b,c,d,e\}$, this can be replaced with four pairs $ab,bc,cd$ and $da$ and two triples 
$ace$ and $bde$. There are fifteen distinct ways to do this, and the resulting systems may or may not be isomorphic. 
Conversely, an MMPTS($v$) of Type-Q with leave $(a,b,c,d)$, and triples $ace$ and $bde$ can be converted to a 
PBD($v,\{3,5^*\}$) having precisely one block of size 5. We take advantage of this correspondence in Section \ref{PBD17} 
to enumerate PBD($17,\{3,5^*\}$)s.

In the analysis of MMPTS(17)s (and other designs), the numbers of Pasch and mitre configurations are determined. A 
\emph{Pasch} configuration in a PTS($v$), also known as a \emph{quadrilateral}, comprises a set of four triples of the 
system on six points, having the form $abc, ade, bdf,cef$. A \emph{mitre} configuration has five triples on seven 
points, having the form $abc, ade, afg, bdf, ceg$; the point $a$ that appears in three triples will be called the 
\emph{root} of the mitre and the other six points will be called the \emph{leafs}. An $(r,s)$-configuration in a 
PTS($v$) is a configuration having $r$ triples and $s$ points. So a Pasch configuration is a $(4,6)$-configuration, and 
in fact the only one. A mitre is a $(5,7)$-configuration and there is one other, obtained from a Pasch configuration by 
adding a triple $afg$. It is a longstanding conjecture regarding Steiner triple systems that, given $r\ge 4$, there 
exists $v_r$ such that for all admissible $v\ge v_r$, there exists an STS($v$) having no $(r',r'+2)$-configurations for 
all $r'$ with $4\le r' \le r$ \cite{E}. Such an STS($v$) is said to be \emph{$r$-sparse}. A 4-sparse system is therefore 
one without Pasch configurations, and is also known as an \emph{anti-Pasch} or \emph{quadrilateral-free} system. An 
\emph{anti-mitre} system is one without mitres, and a 5-sparse system is one without both mitres and Pasch 
configurations, equivalently it is both anti-Pasch and anti-mitre. It is known that there are anti-Pasch STS($v$)s for 
every admissible $v\ne 7,13$ \cite{LCGG,GGW}, and an anti-mitre STS($v$) for every admissible $v\ne 9$ \cite{W2}. It 
was shown in \cite{W1} that 5-sparse STS($v$) exist for almost all admissible $v$, and their existence for orders 
$v\equiv 3$ (mod 6) ($v\ne 9,15$) was established in \cite{W3}. Further results are given in \cite{F}, but the exact 
set of values $v$ for which 5-sparse STS($v$)s exist remains undetermined.

Given an MMPTS(17) and one of the points of its leave, by deleting this point and the 7 triples that contain it, an 
MMPTS(16) is formed. Since all MMPTS(16)s may be formed in this way, it is possible to determine all such systems 
from the MMPTS(17)s. After applying an isomorphism check, we establish that there are 47\,744\,568 MMPTS(16)s. Details 
of these systems are given in Section \ref{Determining16}.

\section{Searching for MMPTS(17)s}\label{Searching17}
\subsection{Seeding the search}\label{Seeding}
The point set of each design was taken to be $V=\{0,1,\ldots,16\}$ and the leave was taken to be $\{01,03,12,23\}$.  
With equipment available to us at the present time it was not feasible to search all possible triples on $V$ 
exhaustively, and this would anyway have been a colossal waste of resources. Instead, a small number of partial 
designs (\emph{seeds}) were constructed that, without loss of generality, could potentially be extended to provide 
collectively at least one representative of each isomorphism class of MMPTS(17)s. These seeds were formed as follows.

Consider first the non-quintuple case. Without loss of generality, the leave and some of the triples of a non-quintuple 
MMPTS(17) may be taken as:
\[\{0,1\},\{0,3\},\{1,2\},\{2,3\},\{0,2,4\},\{1,3,5\},\{0,5,6\},\] 
\[\{0,7,8\},\{0,9,10\},\{0,11,12\},\{0,13,14\},\{0,15,16\}, \{2,5,7\}.\]
Again without loss of generality, the design may be taken to contain either the triple $\{2,6,8\}$ or $\{2,6,9\}$. For 
the first of these possibilities, the next triple may be taken to be $\{2,9,11\}$, while for the second possibility 
the next triple may be taken as either $\{2,8,10\}$ or $\{2,8,11\}$, all without loss of generality. The design 
requires three further blocks containing the point 2. By considering the possibilities it will be seen that, without 
loss of generality, the design must be an extension of one of the five possibilities 
$N_1,N_2,N_3,N_4$ and $N_5$ shown in Table \ref{Table 2}.

Next consider the quintuple case. The argument is almost identical to the non-quintuple case. The only difference is 
that the triple $\{1,3,5\}$ is replaced by the triple $\{1,3,4\}$. Without loss of generality, the design must be an 
extension of one of the five possibilities $Q_1,Q_2,Q_3,Q_{2a}$ and $Q_4$ shown in Table \ref{Table 3}. The 
reason that the fourth of these alternatives is denoted $Q_{2a}$ is that there exists an isomorphism of this design and 
$Q_2$. An isomorphism taking $Q_2$ to $Q_{2a}$ is given by the permutation $(5,13)(6,14)(7,15,10)(8,16,9)$. 
Consequently, design $Q_{2a}$ is not required as a seed in an exhaustive search.

\begin{table}[htbp]
\begin{center}
 \[\begin{array}{|c|c|c|c|c|}\hline N_1 & N_2 & N_3 & N_4 & N_5\\ \hline 0,2,4 & 0,2,4 & 0,2,4 & 0,2,4 & 0,2,4\\ 
1,3,5 & 1,3,5 & 1,3,5 & 1,3,5 & 1,3,5\\ 0,5,6 & 0,5,6 & 0,5,6 & 0,5,6 & 0,5,6\\ 
0,7,8 & 0,7,8 & 0,7,8 & 0,7,8 & 0,7,8\\ 0,9,10 & 0,9,10 & 0,9,10 & 0,9,10 & 0,9,10\\ 
0,11,12 & 0,11,12 & 0,11,12 & 0,11,12 & 0,11,12\\ 0,13,14 & 0,13,14 & 0,13,14\ & 0,13,14\ & 0,13,14\\\ 
0,15,16 & 0,15,16 & 0,15,16 & 0,15,16 & 0,15,16\\ 2,5,7 & 2,5,7 & 2,5,7 & 2,5,7 & 2,5,7\\ 
2,6,8 & 2,6,8 & 2,6,9 & 2,6,9 & 2,6,9\\ 2,9,11 & 2,9,11 & 2,8,10 & 2,8,11 & 2,8,11\\ 
2,10,12 & 2,10,13 & 2,11,13 & 2,10,12 & 2,10,13\\ 2,13,15 & 2,12,15 & 2,12,15 & 2,13,15 & 2,12,15\\ 
2,14,16 & 2,14,16 & 2,14,16 & 2,14,16 & 2,14,16\\ \hline
\end{array}\]
\caption{Seeds for the computer search in the non-quintuple case.\label{Table 2}}
\end{center}
\end{table}

\begin{table}[htbp]
\begin{center}
 \[\begin{array}{|c|c|c|c|c|}\hline Q_1 & Q_2 & Q_3 & Q_{2a} & Q_4\\ \hline  0,2,4 & 0,2,4 & 0,2,4 & 0,2,4 & 0,2,4\\ 
1,3,4 & 1,3,4 & 1,3,4 & 1,3,4 & 1,3,4\\ 0,5,6 & 0,5,6 & 0,5,6 & 0,5,6 & 0,5,6\\ 
0,7,8 & 0,7,8 & 0,7,8 & 0,7,8 & 0,7,8\\ 0,9,10 & 0,9,10 & 0,9,10 & 0,9,10 & 0,9,10\\ 
0,11,12 & 0,11,12 & 0,11,12 & 0,11,12 & 0,11,12\\ 0,13,14 & 0,13,14 & 0,13,14\ & 0,13,14\ & 0,13,14\\\ 
0,15,16 & 0,15,16 & 0,15,16 & 0,15,16 & 0,15,16\\ 2,5,7 & 2,5,7 & 2,5,7 & 2,5,7 & 2,5,7\\ 
2,6,8 & 2,6,8 & 2,6,9 & 2,6,9 & 2,6,9\\ 2,9,11 & 2,9,11 & 2,8,10 & 2,8,11 & 2,8,11\\ 
2,10,12 & 2,10,13 & 2,11,13 & 2,10,12 & 2,10,13\\ 2,13,15 & 2,12,15 & 2,12,15 & 2,13,15 & 2,12,15\\ 
2,14,16 & 2,14,16 & 2,14,16 & 2,14,16 & 2,14,16\\ \hline
\end{array}\]
\caption{Seeds for the computer search in the quintuple case.\label{Table 3}}
\end{center}
\end{table}

\subsection{The searching programs}\label{programs}
Given one of the nine seeds described in the previous section, the problem is to cover all uncovered pairs (apart from 
those in the leave) by triples. This can be represented as a graph covering problem by taking a graph $G$ on the points 
$0,1,\ldots,16$ whose edges are the pairs that are neither in the leave nor covered by triples in the seed. A solution 
corresponds to covering each edge of $G$ precisely once with triangles of edges from $G$, and the program must find 
all possible solutions. A program written by the current authors was used for this purpose. If a point is reached where 
an uncovered pair has no remaining triples available to cover it, the program backtracks. If it reaches a point where 
all pairs are covered, it records the solution and then backtracks. The principal difficulty is adjusting the records 
for backtracks. In each of the nine cases, the program took about 90 minutes on a desktop computer to enumerate the 
solutions. It was then adapted to test solutions for isomorphisms and to record solutions as described below. It is 
prudent to verify computational results independently and so a minor variant of Knuth's \emph{dancing links} 
graph covering program \cite{Knuth1, Knuth2} was used to check the number of solutions from each of the 9 runs; each of 
these checks took about 8 minutes. 

The outcomes of the nine runs are shown in Table \ref{Table 4}.
\begin{table}[htbp]
\begin{center}
\[\begin{array}[t]{|l|c|}\hline \mbox{Seed} & \mbox{Number of designs} \\ \hline 
   N_1 & 161\,885\,696 \\ N_2 & 165\,881\,472 \\ N_3 & 166\,118\,112 \\ N_4 & 170\,428\,416 \\ N_5 & 171\,571\,376 \\
   \hline \mbox{Total} & 835\,885\,072\\ \hline  \end{array} ~~~~~
  \begin{array}[t]{|l|c|}\hline \mbox{Seed} & \mbox{Number of designs} \\ \hline 
   Q_1 & 134\,263\,296 \\ Q_2 & 140\,978\,304 \\ Q_3 & 142\,761\,312 \\ Q_4 & 144\,371\,376 \\  \hline \mbox{Total} & 
   562\,374\,288\\ \hline  \end{array}  \]
\caption{Search results for the numbers of MMPTS(17)s.\label{Table 4}}
\end{center}
\end{table}
The total number of systems obtained in the nine runs was $1\,398\,259\,360$. No solution for the quintuple case 
can be isomorphic with any solution from the non-quintuple case, but within each class there are many isomorphisms. 

\subsection{Isomorphism testing}\label{testing17}
The software package \emph{nauty} \cite{M} was used to analyse the search outputs for isomorphisms. This treats each 
solution as a bipartite graph obtained from the design with an edge joining point $x$ to triple $B$ if and only if 
$x\in B$. Two designs are isomorphic if and only if the two bipartite graphs are isomorphic. Nauty assigns a canonical 
labelling to each graph so that two graphs are isomorphic if and only if they have the same canonical labelling. The 
canonical labellings were stored, as they were produced, in a sorted list. The time taken to insert a new labelling 
into such a list is proportional to the logarithm of the number of objects in the list. If the labelling was already 
present in the list, the system was immediately recognised as an isomorphic copy and so could be discarded. Nauty is a 
long-established package with an excellent record for reliability, so there is a high degree of confidence in the 
results. Nevertheless, steps were taken to check these results. 

Given a solution produced by the searching program, three numbers were obtained for each point $x$ of the 
design: the number $n_p(x)$ of Pasch configurations incident with $x$, the number of mitres $n_r(x)$ with $x$ as a 
root, and the number of mitres $n_l(x)$ with $x$ as a leaf. In the quintuple case, any isomorphism between two solutions 
must preserve the point sets $V_1=\{0,1,2,3\},V_2=\{4\},V_3=\{5,6\ldots,16\}$, and in the non-quintuple case any 
isomorphism must preserve the point sets $V_1=\{0,1,2,3\},V_2=\{4,5\}, V_3=\{6,7\ldots,16\}$. Consequently the three 
sets of ordered triples $\{(n_p(x),n_r(x),n_l(x)): x\in V_i\}$ for $i=1,2,3$ form an invariant in each case. Using this 
invariant, it was easy to determine whether or not two solutions with the same invariant were isomorphic since the 
number of possible isomorphisms was small in every case. Samples of the search output were taken and the results 
produced by consideration of this invariant were found to agree with those produced by nauty. It had been hoped that 
this invariant might be complete, but this is not the case: there are non-isomorphic solutions sharing the same 
invariant.

Table \ref{Table 5} shows the number of isomorphism classes in the two cases. Altogether, up to isomorphism, there are
precisely 35\,810\,097 MMPTS(17)s.
\begin{table}[htbp]
\begin{center}
\[\begin{array}{|l|r|} \hline \mbox{Non-quintuple case:} & 33\,459\,364 \\ \hline
 \mbox{Quintuple case:} & 2\,350\,733 \\ \hline
 \mbox{Total number:} & 35\,810\,097 \\ \hline
\end{array}\]
\caption{The numbers of isomorphism classes.\label{Table 5}}
\end{center}
\end{table}

\subsection{Storing the solutions}\label{Storing17}
Storing a representative of each isomorphism class in the form of 44 triples requires about 14 GB of storage space. 
However this can be brought down to a modest 1.5 GB as follows. Given one of these designs, each triple 
$\{i,j,k\}$ is ordered so that $i<j<k$. These triples are taken in ascending order by their first entry, and then 
by their second. The list of 44 third entries is recorded with numbers represented by characters. The numbers 0-9 
are represented as characters 0-9 and the numbers 10-16 are represented as characters a-g. In this way a design can be 
represented as a character string of length 44. In principle, not all 44 triples are needed because some triples are 
fixed in all 9 runs, but incorporating this saving isn't worth the extra effort. Compressing the resulting files in 
.rar format reduces the required storage space to 13 MB for the quintuple case, and 179 MB for the non-quintuple 
case. The files are available as MMPTS(17)\_Type-Q.rar and MMPTS(17)\_Type-N.rar from \cite{DDG}. 

\section{Properties of the MMPTS(17)s}\label{Properties17}
The number of Pasch configurations, mitres, and STS(7) subsystems (Fano planes) in each design was 
determined by two independently written programs. 

In the quintuple case, 78 designs have no Pasch configurations, and 2 designs have the maximum number of 44 Pasch 
configurations. In the non-quintuple case, 843 designs have no Pasch configurations, and 1 design has the maximum 
number of 47 Pasch configurations. Hence there are anti-Pasch (4-sparse) MMPTS(17)s of both types. Table 
\ref{Table Pasch17} gives the numbers of Type-Q and Type-N MMPTS(17)s having $P$ Pasch configurations for $0\le P\le 
47$. 

\begin{table}[htbp]
  \footnotesize
    \begin{center}
      \begin{tabular}[t]{|c|r|r|}
	\hline
	$P$& Type-Q& Type-N\\ \hline
	0& 78 &843\\
	1& 487 &7\,531\\
	2& 2\,415 &37\,136\\
	3& 8\,161 &124\,920\\
	4& 21\,834 &326\,465\\
	5& 46\,990 &696\,112\\
	6& 85\,332 &1\,258\,629\\
	7& 133\,629 &1\,962\,982\\
	8& 185\,730 &2\,697\,978\\
	9& 227\,426 &3\,292\,317\\
	10& 255\,310 &3\,642\,666\\
	11& 261\,623 &3\,697\,358\\
	12& 248\,943 &3\,477\,394\\
	13& 220\,388 &3\,053\,390\\
	14& 182\,791 &2\,528\,980\\
	15& 143\,166 &1\,984\,888\\
	 \hline  
      \end{tabular}
      \hfill
      \begin{tabular}[t]{|c|r|r|}
	\hline
	$P$& Type-Q& Type-N\\ \hline			
	16& 106\,560 &1\,487\,105\\
	17& 76\,724 &1\,072\,999\\
	18& 52\,071 &745\,422\\
	19& 34\,839 &501\,176\\
	20& 21\,773 &325\,369\\
	21& 13\,707 &206\,993\\
	22& 8\,038 &128\,418\\
	23& 4\,957 &78\,709\\
	24& 2\,782 &48\,182\\
	25& 1\,919 &29\,380\\
	26& 1\,088 &17\,815\\
	27& 676 &11\,220\\
	28& 423 &6\,664\\
	29& 271 &4\,032\\
	30& 208 &2\,456\\
	31& 127 &1\,578\\
	\hline   
      \end{tabular}
      \hfill
      \begin{tabular}[t]{|c|r|r|}
	\hline
	$P$& Type-Q& Type-N\\ \hline
	32& 105 &868\\
	33& 35 &578\\
	34& 39 &303\\
	35& 42 &177\\
	36& 9 &133\\
	37& 4 &87\\
	38& 11 &32\\
	39& 9 &34\\
	40& 6 &6\\
	41& 1 &17\\
	42& 0 &14\\
	43& 4 &1\\
	44& 2 &1\\
	45& 0 &4\\
	46& 0 &1\\
	47& 0 &1\\
	\hline   
      \end{tabular}
    \caption{Pasch configuration frequencies in MMPTS(17)s\label{Table Pasch17}}
  \end{center}
\end{table}

\vspace{-5mm}
In the quintuple case, 4 designs have the minimum number of 5 mitres, and 3 designs have the maximum number of 48 
mitres. In the non-quintuple case, 1 design has the minimum number of 3 mitres, and 1 design has the maximum number of 
47 mitres. Note there are no anti-mitre and, consequently, no 5-sparse MMPTS(17)s.  Table \ref{Table Mitre17} gives the 
numbers of Type-Q and Type-N MMPTS(17)s having $M$ mitres for $0\le M\le 48$.

\begin{table}[htbp]
  \footnotesize
    \begin{center}
      \begin{tabular}[t]{|c|r|r|}
	\hline
	$M$& Type-Q& Type-N\\ \hline
	0& 0 &0\\
	1& 0 &0\\
	2& 0 &0\\
	3& 0 &1\\
	4& 0 &1\\
	5& 4 &13\\
	6& 13 &55\\
	7& 21 &195\\
	8& 69 &695\\
	9& 224 &2\,368\\
	10& 533 &6\,855\\
	11& 1\,367 &18\,010\\
	12& 3\,406 &43\,678\\
	13& 7\,139 &98\,058\\
	14& 14\,637 &199\,350\\
	15&  26\,544 &373\,500\\	
	\hline   
      \end{tabular}
      \hfill
      \begin{tabular}[t]{|c|r|r|}
	\hline
	$M$& Type-Q& Type-N\\ \hline			
	16&44\,626  &643\,122\\
	17& 69\,604 &1\,020\,033\\
	18& 100\,089  &1\,493\,054\\
	19& 134\,564 &2\,026\,513\\
	20& 167\,649 &2\,560\,458\\
	21& 197\,183 &3\,014\,058\\
	22& 217\,341 &3\,317\,425\\
	23& 225\,116 &3\,422\,818\\
	24& 220\,267 &3\,300\,754\\
	25& 205\,106 &2\,992\,985\\
	26& 180\,070 &2\,547\,896\\
	27& 150\,370 &2\,042\,250\\
	28& 119\,906 &1\,536\,133\\
	29& 89\,588 &1\,086\,488\\
	30& 64\,205 &720\,441\\
	31& 43\,397 &449\,377\\
	\hline   
      \end{tabular}	
      \hfill
      \begin{tabular}[t]{|c|r|r|}
	\hline
	$M$& Type-Q& Type-N\\ \hline	
	32& 27\,912 &263\,661\\
	33& 17\,221 &143\,731\\	
	34& 10\,441 &74\,012\\
	35&5\,749  &35\,097\\
	36& 3\,131 &15\,668\\
	37& 1\,640 &6\,630\\
	38& 852 &2\,498\\
	39& 376 &985\\
	40& 201 &320\\
	41& 73 &101\\
	42& 56 &46\\
	43& 12 &19\\
	44& 24 &10\\
	45& 3 &1\\
	46& 1 &0\\
	47& 0 &1\\
	48& 3 &0\\
	\hline   
      \end{tabular}
    \caption{Mitre frequencies in MMPTS(17)s \label{Table Mitre17}}
  \end{center}
\end{table}

Table \ref{Table Fano17} gives the numbers of Type-Q and Type-N MMPTS(17)s having $S$ STS(7) subsystems (also known as 
Fano planes) for $0\le S\le 4$ (the maximum number in any MMPTS(17)).
\begin{table}[htbp]
\begin{center}
\begin{tabular}{|c|r|r|} \hline $S$  & Type-Q & Type-N\\ \hline
  0 & 2\,335\,059 & 33\,198\,472 \\
  1 & 15\,499 & 258\,750 \\
  2 & 125 & 2\,104 \\
  3 & 46 & 38 \\
  4 & 4 & 0 \\ \hline
\end{tabular}
\caption{STS(7) frequencies in MMPTS(17)s.\label{Table Fano17}}
\end{center}
\end{table}

The automorphism groups of the designs were also determined. The numbers of each order were checked by two 
independently written programs. Table \ref{Table Auto17} gives full details. The results in Table \ref{Table Auto17} 
can be also used to determine the total number of MMPTS(17)s on a specified point set from the orbit-stabiliser theorem 
as
\[\sum_D \frac{17!}{|Aut(D)|}, \] 
where $Aut(D)$ denotes the automorphism group of the design $D$, and $D$ runs through the representatives of the 
isomorphism classes of the MMPTS(17)s. The number of Type-Q is 
834\,485\,388\,804\,587\,520\,000 and the number of Type-N is 11\,900\,365\,174\,781\,411\,328\,000. Hence the total 
number of MMPTS(17)s is 12\,734\,850\,563\,585\,998\,848\,000.

\begin{table}[htbp]
\begin{center}
\begin{tabular}{|c|c|r|r|} \hline \multicolumn{2}{|c|}{Automorphisms} & Type-Q & Type-N\\ 
\cline{1-2}  order & group & number of systems & number of systems \\ \hline
 12 & dihedral $D_{12}$ & 3 & 0 \\
 8 & dihedral $D_8$ & 20 & 0 \\
 8 & $\mathbb{Z}_4\times \mathbb{Z}_2$ & 4 & 0 \\
 6 & dihedral $D_6$ ($=S_3$) & 27 & 0 \\
 6 & $\mathbb{Z}_6$ & 3 & 0 \\
 4 & Klein & 269 & 0 \\
 4 & $\mathbb{Z}_4$ & 74 & 0 \\
 3 & $\mathbb{Z}_3$ & 108 & 0 \\
 2 & $\mathbf{Z}_2$ & 8\,470 & 3\,992 \\
 1 & trivial & 2\,341\,755 & 33\,455\,372 \\ \hline
\end{tabular}
\caption{Automorphisms of MMPTS(17)s.\label{Table Auto17}}
\end{center}
\end{table}

\section{The PBD$(17,\{3,5^*\})$s}\label{PBD17}
The PBD$(17,\{3,5^*\})$s were enumerated by two methods. The first was to take the 2\,350\,733 non-isomorphic 
MMPTS(17)s of Type-Q, delete the blocks 024 and 134 (giving a leave forming the complete graph $K_5$ on the set 
$\{0,1,2,3,4\}$), and then test the resulting systems for isomorphisms. The second was to start again from Table 
\ref{Table 3}, delete the blocks 024, 134 from each seed, generate systems from the seeds, and then test for 
isomorphisms. Isomorphisms were identified using nauty and checked with our own independent program. Both approaches 
gave the number of non-isomorphic PBD$(17,\{3,5^*\})$s as 157\,151.  Thus the vast majority of these designs must each 
generate 15 non-isomorphic MMPTS(17)s when two intersecting triples are formed from the leave $K_5$ in all 
possible ways. The non-isomorphic PBD$(17,\{3,5^*\})$s are available in compact form as strings of length 42 from 
\cite{DDG} in the file PBD(17,$\{$3,5$\}$).rar (1 MB).

Analysis of these systems shows that there are 39 systems with no Pasch configurations, and 1 system with the maximum 
number of 34 Pasch configurations. Table \ref{Table 17PBDPasch} gives the number $n_p$ of systems having $P$ Pasch 
configurations for $0\le P\le 34$.

\begin{table}[htbp]
\footnotesize
  \begin{center}
    \begin{tabular}[t]{|c|r|}
      \hline
      $P$& $n_p$\\ \hline
      0&39\\
      1&259\\
      2&911\\
      3&2\,701\\
      4&5\,440\\
      5&9\,324\\
      6&13\,619\\
      7&17\,149\\
      8&19\,363\\
      \hline  
    \end{tabular}
    \hfill
    \begin{tabular}[t]{|c|r|}
      \hline
      $P$& $n_p$\\ \hline 
      9&19\,285\\
      10&17\,904\\
      11&14\,844\\
      12&11\,935\\
      13&8\,479\\
      14&5\,958\\
      15&3\,982\\
      16&2\,468\\
      17&1\,467\\
      \hline   
    \end{tabular}
    \hfill
    \begin{tabular}[t]{|c|r|}
      \hline
      $P$& $n_p$\\ \hline
      18&855\\
      19&472\\
      20&270\\
      21&176\\
      22&88\\
      23&61\\
      24&35\\
      25&23\\
      26&9\\
      \hline  
    \end{tabular}
    \hfill
    \begin{tabular}[t]{|c|r|}
      \hline
      $P$& $n_p$\\ \hline
      27&8\\
      28&12\\
      29&1\\
      30&5\\
      31&5\\
      32&1\\
      33&2\\
      34&1\\
      \hline  
    \end{tabular}
    \caption{Pasch configuration frequencies in PBD$(17,\{3,5^*\})$s.\label{Table 17PBDPasch}}
  \end{center}
\end{table}

The minimum number of mitres in a PBD$(17,\{3,5^*\})$ is 5 (10 systems have this number), and 2 systems have the 
maximum number of 36 mitres. Table \ref{Table 17PBDMitre} gives the number $n_m$ of systems having $M$ mitres for $5\le 
M\le 36$. 

\begin{table}[htbp]
\footnotesize
  \begin{center}
    \begin{tabular}[t]{|c|r|}
      \hline
      $M$& $n_m$\\ \hline
      5&10\\
      6&19\\
      7&52\\
      8&153\\
      9&430\\
      10&1\,001\\
      11&2\,211\\
      12&4\,026\\
      \hline   
    \end{tabular}
    \hfill
    \begin{tabular}[t]{|c|r|}
      \hline
      $M$& $n_m$\\ \hline
      13&6\,640\\
      14&10\,146\\
      15&13\,561\\
      16&16\,615\\
      17&18\,614\\
      18&18\,821\\
      19&17\,540\\
      20&15\,142\\
      \hline 
    \end{tabular}
    \hfill
    \begin{tabular}[t]{|c|r|}
      \hline
      $M$& $n_m$\\ \hline
      21&11\,618\\
      22&8\,449\\
      23&5\,448\\
      24&3\,306\\
      25&1\,809\\
      26&892\\
      27&386\\
      28&171\\
      \hline  
    \end{tabular}
    \hfill
    \begin{tabular}[t]{|c|r|}
      \hline
      $M$& $n_m$\\ \hline
      29&54\\
      30&22\\
      31&2\\
      32&9\\
      34&1\\
      35&1\\
      36&2\\
      \hline  
    \end{tabular}
    \caption{Mitre frequencies in PBD$(17,\{3,5^*\})$s. \label{Table 17PBDMitre}}
  \end{center}
\end{table}

Each PBD$(17,\{3,5^*\})$ has either $S=0$ or $S=1$ STS(7) subsystems (Fano planes). Table \ref{Table STS7in17PBD} gives 
the number $n_s$ of systems having $S$ STS(7) subsystems.

\begin{table}[htbp]
  \begin{center}
    \begin{tabular}{|c|r|} \hline 
      $S$& $n_s$\\ \hline
      0 & 156\,713 \\
      1 &438 \\\hline
    \end{tabular}
    \caption{STS(7) frequencies in PBD$(17,\{3,5^*\})$s.\label{Table STS7in17PBD}}
  \end{center}
\end{table}

Table \ref{Table 17PBDautom} tabulates the automorphism groups of the PBD$(17,\{3,5^*\})$s. \linebreak Again, the 
numbers of each order were checked by two independently written programs. The groups of order 24 were identified using 
the GAP system \cite{GAP}. The symbol $\rtimes$ denotes a semidirect product. The total number of PBD$(17,\{3,5^*\})$s  
on a specified point set can be obtained from the orbit-stabiliser theorem, and is 55\,632\,359\,253\,639\,168\,000.

\begin{table}[htbp]
  \begin{center}
    \begin{tabular}{|c|c|r|} 
      \hline \multicolumn{2}{|c|}{Automorphisms} & number of PBD$(17,\{3,5^*\})$s \\ 
      \cline{1-2} order & group & \\ \hline
      24 & $S_4$ & 2 \\
      24 & $\mathbb{Z}_4\times S_3$ & 1 \\
      24 & $(\mathbb{Z}_6\times\mathbb{Z}_2)\rtimes\mathbb{Z}_2$ & 1 \\
      12 & $\mathbb{Z}_{12}$ & 2 \\
      12 & dihedral $D_{12}$ & 2 \\
      12 & $\mathbb{Z}_3\rtimes\mathbb{Z}_4$ & 2 \\
      8 & dihedral $D_8$ & 5 \\
      8 & $\mathbb{Z}_4\times\mathbb{Z}_2$ & 3 \\
      6 & $\mathbb{Z}_6$ & 15 \\
      6 & dihedral $D_6~(=S_3)$ & 18 \\
      4 & $\mathbb{Z}_4$ & 19 \\
      4 & Klein & 63 \\
      3 & $\mathbf{Z}_3$ & 64 \\
      2 & $\mathbf{Z}_2$ & 1\,190 \\
      1 & trivial & 155\,764 \\ \hline
    \end{tabular}
    \caption{Automorphisms of PBD$(17,\{3,5^*\})$s.\label{Table 17PBDautom}}
  \end{center}
\end{table}

\section{The MMPTS(16)s}\label{Determining16}
As explained in the Introduction, the MMPTS(16)s can be obtained by deleting one of the points 0, 1, 2 or 3 of each of 
the MMPTS(17)s described in Section \ref{Searching17}. This results in $4\times 35\,810\,097=143\,240\,388$ systems. 
But amongst these there are isomorphisms, and nauty was used to eliminate isomorphic copies. It was found that there 
are 47\,744\,568 non-isomorphic solutions. The dancing links program was used to check this number with an independent 
search. Of these solutions, 3\,133\,686 originate from MMPTS(17)s of Type-Q, and the remaining 44\,610\,882 
are from systems of Type-N. If an MMPTS(16) originates from an MMPTS(17) of Type-Q, then adding a new point to 
each of the seven non-intersecting pairs and to one of the three intersecting pairs of its leave gives a Type-Q 
MMPTS(17). Hence no MMPTS(16) originating from a quintuple type MMPTS(17) can be isomorphic to one originating from a 
non-quintuple type MMPTS(17). A MMPTS(16) derived from a quintuple type MMPTS(17) will be described as Type-Q, and a 
MMPTS(16) derived from a non-quintuple type MMPTS(17) will be described as Type-N. The designs are available in files
MMPTS(16)\_Type-Q.rar (17 MB) and MMPTS(16)\_Type-N.rar (208 MB) from \cite{DDG}. The designs are mapped so 
that the leave is $\{01,02,03,45,67,\ldots,\textrm{ef}\}$, with characters a to f representing points 10 to 15, and the 
designs are specified in compact form as strings of length 37.

For MMPTS(16)s of Type-Q, 1\,346 designs have no Pasch configurations, and 1 design has the maximum number of 43 Pasch 
configurations. For MMPTS(16)s of type-N, 15\,877 designs have no Pasch configurations, and 1 design has the maximum 
number of 44 Pasch configurations. Hence there are anti-Pasch (4-sparse) MMPTS(16)s of both types. Table \ref{Table 
Pasch16} gives the numbers of Type-Q and Type-N MMPTS(16)s having $P$ Pasch configurations for $0\le P\le 44$.

\begin{table}[htbp]
  \footnotesize
    \begin{center}
      \begin{tabular}[t]{|c|r|r|}
	\hline
	$P$& Type-Q& Type-N\\ \hline
	0&1\,346&15\,877\\
	1&9\,571&125\,270\\
	2&38\,020&502\,936\\
	3&100\,154&1\,347\,901\\
	4&198\,137&2\,704\,815\\
	5&309\,441&4\,278\,201\\
	6&395\,447&5\,549\,182\\
	7&432\,141&6\,121\,166\\
	8&414\,833&5\,905\,615\\
	9&358\,284&5\,122\,025\\
	10&282\,328&4\,069\,661\\
	11&206\,951&3\,016\,887\\
	12&142\,368&2\,113\,200\\
	13&94\,556&1\,407\,703\\
	14&59\,617&904\,747\\
	\hline   
      \end{tabular}
      \hfill
      \begin{tabular}[t]{|c|r|r|}
	\hline
	$P$& Type-Q& Type-N\\ \hline			
	15&36\,534&565\,757\\
	16&22\,010&344\,867\\
	17&12\,897&207\,237\\
	18&7\,562&123\,640\\
	19&4\,448&73\,743\\
	20&2\,606&43\,954\\
	21&1\,661&26\,351\\
	22&1\,013&15\,917\\
	23&653&9\,497\\
	24&355&5\,761\\
	25&264&3\,480\\
	26&148&2\,096\\
	27&120&1\,327\\
	28&85&755\\
	29&49&510\\
	\hline   
      \end{tabular}
      \hfill
      \begin{tabular}[t]{|c|r|r|}
	\hline
	$P$& Type-Q& Type-N\\ \hline
	30&29&314\\
	31&18&163\\
	32&14&125\\
	33&10&76\\
	34&2&38\\
	35&5&36\\
	36&1&17\\
	37&1&12\\
	38&2&10\\
	39&2&3\\
	40&1&3\\
	41&1&5\\
	42&0&1\\
	43&1&0\\
	44&0&1\\
	\hline   
      \end{tabular}
    \caption{Pasch configuration frequencies in MMPTS(16)s \label{Table Pasch16}}
  \end{center}
\end{table}
For MMPTS(16)s of Type-Q, 2 designs have no mitres, and 1 design has the maximum number of 48 mitres. For MMPTS(16)s of 
type-N, 5 designs have no mitres, and 1 design has the maximum number of 45 mitres. Hence there are anti-mitre 
MMPTS(16)s of both types. Table \ref{Table Mitre16} gives the numbers of Type-Q and Type-N MMPTS(16)s having $M$ 
mitres for $0\le M\le 48$. The 7 designs with no mitres all have Pasch configurations, so there are no 5-sparse 
MMPTS(16)s. 

\begin{table}[htbp]
  \footnotesize
    \begin{center}
      \begin{tabular}[t]{|c|r|r|}
	\hline
	$M$& Type-Q& Type-N\\ \hline
	0&2&5\\
	1&3&14\\
	2&26&119\\
	3&116&879\\
	4&526&5\,005\\
	5&2\,165&21\,590\\
	6&6\,559&73\,898\\
	7&17\,563&206\,848\\
	8&39\,376&488\,364\\
	9&77\,123&988\,266\\
	10&131\,409&1\,741\,921\\
	11&198\,113&2\,700\,157\\
	12&267\,765&3\,725\,024\\
	13&324\,803&4\,620\,691\\
	14&359\,205&5\,184\,315\\
	15&361\,423&5\,303\,746\\
	\hline   
      \end{tabular}
      \hfill
      \begin{tabular}[t]{|c|r|r|}
	\hline
	$M$& Type-Q& Type-N\\ \hline			
	16&335\,963&4\,967\,389\\
	17&289\,977&4\,289\,496\\
	18&231\,425&3\,427\,280\\
	19&172\,907&2\,532\,692\\
	20&121\,074&1\,746\,307\\
	21&79\,603&1\,122\,077\\
	22&50\,082&677\,356\\
	23&29\,615&385\,028\\
	24&16\,792&205\,365\\
	25&9\,284&103\,465\\
	26&4\,974&50\,430\\
	27&2\,591&23\,424\\
	28&1\,548&10\,649\\
	29&695&4\,811\\
	30&441&2\,148\\
	31&195&936\\
	\hline   
      \end{tabular}
      \hfill
      \begin{tabular}[t]{|c|r|r|}
	\hline
	$M$& Type-Q& Type-N\\ \hline
	32&153&469\\	
	33&41&310\\
	34&77&153\\
	35&4&89\\
	36&35&87\\
	37&8&16\\
	38&10&23\\
	39&1&15\\
	40&9&3\\
	41&0&6\\
	42&3&14\\
	43&0&1\\
	44&1&0\\
	45&0&1\\
	46&0&0\\
	47&0&0\\
	48&1&0\\
	\hline   
      \end{tabular}
    \caption{Mitre frequencies in MMPTS(16)s \label{Table Mitre16}}
  \end{center}
\end{table}

Table \ref{Table Fano16} gives the numbers of Type-Q and Type-N MMPTS(16)s having $S$ STS(7) subsystems (Fano planes) 
for $0\le S\le 4$ (the maximum number in any MMPTS(16)).
\begin{table}[htbp]
\begin{center}
\begin{tabular}{|c|r|r|} \hline $S$  & Type-Q & Type-N\\ \hline
  0 & 3\,120\,879 & 44\,379\,670 \\
  1 & 12\,703 & 229\,796 \\
  2 & 80 & 1\,392 \\
  3 & 22 & 24 \\
  4 & 2 & 0 \\ \hline
\end{tabular}
\caption{STS(7) frequencies in MMPTS(16)s.\label{Table Fano16}}
\end{center}
\end{table}

Table \ref{Table Auto16} tabulates the automorphism groups of the MMPTS(16)s. Again, the numbers of each order were 
checked by two independently written programs. The total number of MMPTS(16)s on a specified point set can be obtained 
from the orbit-stabiliser theorem, and this gives the number of Type-Q systems as 65\,449\,834\,416\,046\,080\,000 and 
the number of Type-N systems as 933\,361\,974\,492\,659\,712\,000. Thus the overall total number of MMPTS(16)s is 
998\,811\,808\,908\,705\,792\,000.

\begin{table}[htbp]
\begin{center}
\begin{tabular}{|c|c|r|r|} \hline \multicolumn{2}{|c|}{Automorphisms} & Type-Q & Type-N\\ 
\cline{1-2}  order & group & number of systems & number of systems \\ \hline
 12 & dihedral $D_{12}$ & 9 & 0 \\
 12 & alternating $A_4$ & 2 & 0 \\
 12 & $\mathbb{Z}_6\times \mathbb{Z}_2$ & 1 & 0 \\
 6 & dihedral $D_6$ ($=S_3$) & 33 & 6 \\
 6 & $\mathbb{Z}_6$ & 28 & 0 \\
 4 & Klein & 470 & 0 \\
 4 & $\mathbb{Z}_4$ & 32 & 0 \\
 3 & $\mathbb{Z}_3$ & 280 & 504 \\
 2 & $\mathbf{Z}_2$ & 9\,802 & 1\,434 \\
 1 & trivial & 3\,123\,029 & 44\,608\,938 \\ \hline
\end{tabular}
\caption{Automorphisms of MMPTS(16)s.\label{Table Auto16}}
\end{center}
\end{table}

\section{Concluding remarks}
Details of MMPTS($v$)s for $v\le 11$ are given in \cite{CR1}. In this Section we give some details of the systems of 
orders 12 and 14 that have been enumerated elsewhere. As recorded in \cite{GGQS}, there are precisely five 
non-isomorphic MMPTS(12)s arising from the two non-isomorphic STS(13)s, but these are not explicitly listed. Table 
\ref{Table MMPTS12} gives them here in compact form with the numbers of Pasch configurations and mitres that they 
contain, and their automorphism groups. In each case the leave is $\{01,23,\ldots,\mathrm{ab}\}$ with a, b 
representing points 10, 11. Since no STS(13) has an STS(7) subsystem, neither do these MMPTS(12)s. The total 
number of MMPTS(12)s on a specified point set is 1\,197\,504\,000.
\begin{table}[htbp]
  \begin{center}
    \begin{tabular}{|cccc|}\hline MMPTS(12) & \#Pasch & \#mitre & Automorphism group \\ \hline
      468ab798abb8ab7aa99b   &   4   &   6 & $\mathbb{Z}_2$ \\
      468ab79a6b9baba889ba   &   5   &   3 & $\mathbb{Z}_2$ \\
      468ab78ab96ba9abb89a   &   5   &   3 & dihedral $D_6$ ($=S_3$)\\
      468ab94ba87baa9b8ba9   &   4   &   5 & trivial \\
      468aba498b79bb9ab8aa   &   7   &   0 & $\mathbb{Z}_3$ \\ \hline
    \end{tabular}
    \caption{The non-isomorphic MMPTS(12)s.\label{Table MMPTS12}}
 \end{center}
\end{table}

There are 787 non-isomorphic MMPTS(14)s, and these are enumerated in \cite{Gropp} where they are recorded as 
$(14_6,28_3)$ configurations. Here we give counts for Pasch configurations, mitres and STS(7) subsystems, and details 
of the automorphism groups. The designs are available in the file MMPTS(14).txt (23KB) from \cite{DDG}. They are mapped 
so that the leave is $\{01,23,\ldots,\mathrm{cd}\}$ with characters a to d representing points 10 to 13, and they are 
specified in compact form as strings of length 28.

There are 6 MMPTS(14)s with no Pasch configurations, and 1 MMPTS(14) has the maximum number of 63 Pasch 
configurations. Table \ref{Table 14Pasch} gives the number $n_p$ of MMPTS(14)s having $P$ Pasch configurations for 
$0\le P\le 63$.

\begin{table}
    \begin{center}
      \begin{tabular}[t]{|c|r|}
	\hline
	$P$& $n_p$\\ \hline
	0&6\\
	1&7\\
	2&34\\
	3&94\\
	4&94\\
	5&117\\
	6&72\\
	7&62\\
	8&56\\
	\hline   
      \end{tabular}
      \hspace{5mm}
      \begin{tabular}[t]{|c|r|}
	\hline
	$P$& $n_p$\\ \hline			
	9&36\\
	10&26\\
	11&30\\
	12&13\\
	13&29\\
	14&23\\
	15&9\\
	16&3\\
	17&17\\
	\hline  
      \end{tabular}
      \hspace{5mm}
      \begin{tabular}[t]{|c|r|}
	\hline
	$P$& $n_p$\\ \hline
	19&12\\
	20&7\\
	21&5\\
	22&4\\
	23&10\\
	25&2\\
	26&1\\
	27&2\\
	28&1\\
	\hline
      \end{tabular}
      \hspace{5mm}
      \begin{tabular}[t]{|c|r|}
	\hline
	$P$& $n_p$\\ \hline
	29&2\\
	31&7\\
	33&1\\
	35&1\\
	39&1\\
	43&1\\
	47&1\\
	63&1\\		
	\hline
      \end{tabular}
    \caption{Pasch configuration frequencies in MMPTS(14)s \label{Table 14Pasch}}
  \end{center}
\end{table}

There are 15 MMPTS(14)s with no mitres, and 1 MMPTS(14) has the maximum number of 22 mitres. Table \ref{Table 14Mitre} 
gives the number $n_m$ of MMPTS(14)s having $M$ mitres for $0\le M\le 22$. The 15 designs with no mitres all have Pasch 
configurations, so there are no 5-sparse MMPTS(14)s.

\begin{table}
    \begin{center}
      \begin{tabular}[t]{|c|r|}
	\hline
	$M$& $n_m$\\ \hline
	0&15\\
	2&5\\
	3&8\\
	4&24\\
	5&12\\
	\hline
      \end{tabular}
      \hspace{5mm}
      \begin{tabular}[t]{|c|r|}
	\hline
	$M$& $n_m$\\ \hline	
	6&44\\
	7&35\\		
	8&73\\
	9&87\\
	10&90\\
	\hline
      \end{tabular}
      \hspace{5mm}
      \begin{tabular}[t]{|c|r|}
	\hline
	$M$& $n_m$\\ \hline
	11&96\\
	12&100\\
	13&71\\
	14&51\\
	15&40\\
	\hline   
      \end{tabular}
      \hspace{5mm}
      \begin{tabular}[t]{|c|r|}
	\hline
	$M$& $n_m$\\ \hline
	16&23\\
	17&7\\
	18&1\\
	19&3\\
	20&1\\
	22&1\\
	\hline   
      \end{tabular}
    \caption{Mitre frequencies in MMPTS(14)s \label{Table 14Mitre}}
  \end{center}
\end{table}

Table \ref{Table STS7in14} gives the number $n_s$ of MMPTS(14)s having $S$ STS(7) subsystems (Fano planes) for $0\le 
S\le 8$ (the maximum number in any MMPTS(14)).

\begin{table}[htbp]
  \begin{center}
    \begin{tabular}{|r|r|} 
      \hline 
      $S$ & $n_s$ \\ \hline
      0 & 730 \\
      1 & 43 \\
      2 & 11 \\
      4 & 2 \\
      8 & 1 \\ \hline
    \end{tabular}
    \caption{STS(7) frequencies in MMPTS(14)s.\label{Table STS7in14}}
  \end{center}
\end{table}

Table \ref{Table 14autom} tabulates the automorphism groups of the MMPTS(14)s. Again, the numbers of each order were 
checked by two independently written programs. The groups of orders 24 and above were identified using the GAP 
system \cite{GAP}. The symbol $\rtimes$ denotes a semidirect product. The system with the automorphism group of order 
1344 is the system obtained from the projective STS(15) that has an automorphism group $A_8$, the alternating group of 
order 20\,160, and whose automorphism partition consists of a single part containing all 15 points; this is STS(15) 
number 1 in the listing of \cite{MPR}. The total number of MMPTS(14)s on a specified point set can be obtained from the 
orbit-stabiliser theorem, and is 60\,281\,712\,691\,200.

\begin{table}[htbp]
  \begin{center}
    \begin{tabular}{|c|c|r|} \hline 
      \multicolumn{2}{|c|}{Automorphisms} & number of MMPTS(14)s \\  		
      \cline{1-2}  order & group & \\ \hline
	1344 & $\mathbb{Z}_2^3\rtimes PSL(3,2)$ & 1 \\
	192 & $(((\mathbb{Z}_2\times D_8)\rtimes\mathbb{Z}_2)\rtimes\mathbb{Z}_3)\rtimes\mathbb{Z}_2$ & 1 \\
	96 & $(\mathbb{Z}_2^2\times A_4)\rtimes \mathbb{Z}_2$ & 1\\
	32 & $(\mathbb{Z}_2^4)\rtimes \mathbb{Z}_2$ & 3 \\
	24 & symmetric $S_4$ & 3 \\
	24 & $A_4\times\mathbb{Z}_2$ & 2 \\
	21 & $\mathbb{Z}_7\rtimes \mathbb{Z}_3$ & 2 \\
	16 & $\mathbb{Z}_2^4$ & 1 \\
	12 & alternating $A_4$ & 4 \\
	8 & dihedral $D_8$ & 10 \\
	6 & dihedral $D_6$ ($=S_3$) & 1 \\
	6 & $\mathbb{Z}_6$ & 1 \\
	4 & Klein & 18 \\
	4 & $\mathbb{Z}_4$ & 13 \\
	3 & $\mathbf{Z}_3$ & 37 \\
	2 & $\mathbf{Z}_2$ & 40 \\
	1 & trivial & 649 \\ \hline
      \end{tabular}
    \caption{Automorphisms of MMPTS(14)s.\label{Table 14autom}}
  \end{center}
\end{table}

Collecting together previous results and the results of this paper, Table \ref{Table MMPTS(v)} gives the numbers of 
non-isomorphic MMPTS($v$)s for $v\le 19$. The entry for $v=18$ comes from a personal communication \cite{KtoG} cited in 
\cite{Gropp} and given also in \cite[page 354]{HBK}.
\begin{table}[htbp]
  \begin{center}
    \begin{tabular}{|c|c|} 
      \hline $v$  & The number of non-isomorphic MMPTS($v$)s\\ \hline
      $\leq 9$ & 1  \\
      10 & 2  \\
      11 & 2  \\
      12 & 5  \\
      13 & 2  \\
      14 & 787  \\
      15 & 80  \\
      16 & 47\,744\,568  \\
      17 & 35\,810\,097  \\ 
      18 & 210\,611\,385\,743\\
      19 & 11\,084\,874\,829 \\ \hline
    \end{tabular}
  \caption{The number of non-isomorphic MMPTS($v$)s. \label{Table  MMPTS(v)}}
  \end{center}
\end{table}


\newpage

\begin{thebibliography}{99}

\bibitem{HBK} C. J. Colbourn and J. H. Dinitz. The CRC Handbook of Combinatorial Designs (2nd edition), CRC Press,
2007.

\bibitem{CR1} C. J. Colbourn and A. Rosa. Maximal partial Steiner triple systems of order $v\leq 11$. Ars 
Combin. 20, (1985) 5--28.

\bibitem{CR} C. J. Colbourn and A. Rosa. Triple Systems. Oxford University Press, 1999.

\bibitem{CMRS} C. J. Colbourn, E. Mendelsohn, A. Rosa and J. Siran. Anti-mitre Steiner triple systems. Graphs 
Combin. 10 (1994) 215--224.

\bibitem{DDG} F. Demirkale, D. M. Donovan and M. J. Grannell. Maximum partial triple systems on 14, 16 
and 17 points. \texttt{http://www.yildiz.edu.tr/$\sim$fatihd/MMPTS/}

\bibitem{E} P. Erd\H{o}s. Problems and results in combinatorial analysis. Colloquio Internazionale sulle Teorie 
Combinatorie (Rome, 1973), Tomo II, Atti dei Convegni Lincei, No. 17, (1976) 3--17.

\bibitem{F} Y. Fujiwara. Infinite classes of anti-mitre and 5-sparse Steiner triple systems. J. Combin. Des. 14(3) 
(2006) 237-250.

\bibitem{GAP} GAP - Groups, Algorithms, Programming - a System for Computational Discrete Algebra. 
https://www.gap-system.org/

\bibitem{GGQS} M. J. Grannell, T. S. Griggs, K. A. S. Quinn and R. G. Stanton. A census of minimal pair-coverings with 
restricted largest block length. Ars Combin 52 (1999) 71--96.

\bibitem{GGW} M. J. Grannell, T. S. Griggs and C. A. Whitehead. The resolution of the anti-Pasch conjecture. J. 
Combin. Des. 8(4) (2000) 300--309.

\bibitem{Gropp} H. Gropp. Existence and enumeration of configurations. Bayreth. Math. Schr. 74 (2005) 123--129.

\bibitem{KtoG} P. Kaski. Personal communication to Harald Gropp. April 2005.

\bibitem{KO} P. Kaski and P. \"{O}sterg{\aa}rd. The Steiner triple systems of order 19. Math. Comp. 73 (2004) no.248
2075--2092.

\bibitem{Knuth1} D. E. Knuth. Dancing links, \emph{in} Millennial Perspectives in Computer Science: Proceedings of the 
1999 Oxford-Microsoft Symposium in Honour of Sir Tony Hoare. J. Davies, B. Roscoe, and J. Woodcock, editors. Palgrave, 
2000, 187--214.

\bibitem{Knuth2} D. E. Knuth. Selected papers on fun \& games. CSLI Lecture Notes 192, CSLI Publications, Stanford, CA,
2011.

\bibitem{LCGG} A. C. H. Ling, C. J. Colbourn, M. J. Grannell and T. S. Griggs. Construction techniques for anti-Pasch 
Steiner triple systems. J. London Math. Soc. (2), 61 (2000) 641--657.

\bibitem{MPR} R. A. Mathon, K. T. Phelps and A. Rosa. Small Steiner triple systems and their properties. Ars Combin. 
15 (1983) 3--110.

\bibitem{M} B. D. McKay and A. Piperno. Nauty and Traces user’s guide (Version 2.5). Computer Science Department, 
Australian National University, Canberra, Australia (2013).

\bibitem{WCC} H. S. White, F. N. Cole and L. D. Cummings. Complete classification of the triad systems on fifteen 
elements. Memoirs Nat. Acad. Sci. U.S.A. 14 (1919) 1--89.

\bibitem{Wilson} R. M. Wilson. Nonisomorphic Steiner triple systems. Math. Zeitschr. 135 (1973/4) 303--313.

\bibitem{W1} A. J. Wolfe, 5-Sparse Steiner triple systems of order $n$ exist for almost all admissible $n$. Electron. 
J. Combin. 12 (2005) \#R68, 42 pp. (electronic).

\bibitem{W2} A. J. Wolfe. The resolution of the anti-mitre Steiner triple system conjecture. J. Combin. Des. 14(3) 
(2006) 229--236.

\bibitem{W3} A. J. Wolfe, The existence of 5-sparse Steiner triple systems of order $n\equiv 3$ (mod 6), $n\notin\{9, 
15\}$. J. Combin. Theory, Ser. A 115(8) (2008) 1487--1503.


\end{thebibliography}
\end{document}